\documentclass{article}

\font\es=eufm10

\def\gge{\mbox{\es {e}}} 
 
\def\gg{\mbox{\es {g}}}

\def\gC{\mbox{\es {C}}}

\def\gJ{\mbox{\es {J}}}

\def\gP{\mbox{\es {P}}}

\def\so{\mbox{\es {so}}}

\def\spin{\mbox{\es {spin}}}

\def\det{\mbox{\rm {det}}}
\def\diag{\mbox{\rm {diag}}}
\def\Iso{\mbox{\rm {Iso}}}

\def\Ker{\mbox{\rm {Ker}}}
\def\exp{\mbox{\rm {exp}}}
\def\dim{\mbox{\rm {dim}}}
\def\ov{\overline}
\def\wti{\widetilde}
\def\ti{\tilde}

\def\dfrac#1#2{\displaystyle \frac{#1}{#2}}

\def\H{\mbox{\boldmath $H$}}

\def\O{\mbox{\boldmath $O$}}

\def\R{\mbox{\boldmath $R$}}

\def\Z{\mbox{\boldmath $Z$}} 
\def\sR{\mbox{\boldmath $\scriptstyle{R}$}}

\def\0{\mbox{\boldmath {0}}}    
\def\1{\mbox{\boldmath {1}}}      
\def\2{\mbox{\boldmath {2}}}      
\def\3{\mbox{\boldmath {3}}}      
\def\4{\mbox{\boldmath {4}}}      
\def\5{\mbox{\boldmath {5}}}      
\def\6{\mbox{\boldmath {6}}}      
\def\7{\mbox{\boldmath {7}}}      
\def\8{\mbox{\boldmath {8}}}      
\def\9{\mbox{\boldmath {9}}}      
\def\a{\mbox{\boldmath $a$}}

\begin{document}
\setcounter{page}{1} 
\baselineskip=14pt
\begin{center}
{\Large{\bf Fixed points subgroups  by two involutive }}
\end{center}
\begin{center}
{\Large{\bf automorphisms $\sigma$, $\gamma$ of compact exceptional}}
\end{center}
\begin{center}
{\Large{\bf Lie groups $F_4, E_6$ and $E_7$}}
\end{center}
\vspace{2mm}

\begin{center}
{\large{By}}
\end{center}

\begin{center}
{\large{Toshikazu {\sc Miyashita}}}
\end{center}
\vspace{3mm}

{\large{\bf Introduction}}
\vspace{2mm}

For simply connected compact exceptional Lie groups $G = F_4, E_6$ and $E_7$, we consider two involutions $\sigma, \gamma$ and determine the group structure of subgroups $G^{\sigma,\gamma}$  of $G$ which are the intersection $G^\sigma \cap G^{\gamma}$ of the fixed points subgroups of $G^\sigma$ and $G^{\gamma}$. The motivation is as follows. In [1], we determine the group structure of $(F_4)^{\sigma, \sigma'}, (E_6)^{\sigma, \sigma'}$ and $(E_7)^{\sigma, \sigma'}$, and in [2], we also determine the group structure of $(G_2)^{\gamma, \gamma'}, (F_4)^{\gamma, \gamma'}$ and $(E_6)^{\gamma, \gamma'}$. So, in this paper, we try to determine the type of groups $(F_4)^{\sigma, \gamma}, (E_6)^{\sigma, \gamma}$ and $(E_7)^{\sigma, \gamma}$. Our results are the following second columns. The first columns are already known in [3],[4] or [5] and these play an important role to obtain our results. In Table 1, the results of the group structure of $G^{\sigma, \gamma}$ are obtained by the result of $G^\gamma$ and in Table 2, ones are obtained by the result of $G^\sigma$. In this paper, we show the proof of the results of the first and the second line of Table 1 and the third line of Table 2.\\  
 \hspace*{5mm}{\sc Acknowledgment} \quad The author is grateful to Professor Ichiro Yokota for his valuable commentes.\\

\begin{center}
{\bf Table 1}
\end{center}
$$
\begin{array}{lll}
G  &    G^{\gamma}   &   G^{\sigma, \gamma} 
\vspace{2mm} \\ 
\hline
\vspace{-1mm}\\
F_4  &   (Sp(1) \times Sp(3))/\Z_2   &   (Sp(1)\times Sp(1)\times Sp(2))/\Z_2 
\vspace{1mm} \\
E_6  &    (Sp(1)\times SU(6))/\Z_2   &   (Sp(1)\times S(U(2)\times U(4)))/\Z_2 
\vspace{1mm} \\
E_7  &    (SU(2)\times Spin(12))/\Z_2  &  (SU(2)\times Spin(4)\times Spin(8))/(\Z_2 \times \Z_2) 
\end{array}$$

\begin{center}
{\bf Table 2}
\end{center}
$$
\begin{array}{lll}
G  &    G^{\sigma}   &   G^{\sigma, \gamma} 
\vspace{2mm} \\ 
\hline
\vspace{-1mm}\\
F_4  &   Spin(9)   &   (Spin(4)\times Spin(5))/\Z_2 
\vspace{1mm} \\
E_6  &    (U(1)\times Spin(10))/\Z_4   &   (U(1)\times Spin(4)\times Spin(6))/\Z_2 
\vspace{1mm} \\
E_7  &    (SU(2)\times Spin(12))/\Z_2  &  (SU(2)\times Spin(4)\times Spin(8))/(\Z_2 \times \Z_2) 
\end{array}$$

\noindent As for the group $(E_8)^{\sigma, \gamma}$, we can not realize explicitly, however we conjecture 
$$
       (E_8)^{\sigma, \gamma} \cong (Spin(4)\times Spin(12))/(\Z_2\times \Z_2) 
$$
{\sc Remark.} In $E_7$, since $\gamma$ is conjugate to $-\sigma$, we have $(E_7)^\gamma \cong (E_7)^\sigma$. (In detail, see [4].) Note that the results of Table 1 and Table 2 are the same as a set, however they are different as realizations.
\vspace{2mm}

{\large {\bf  Notation}} 
\vspace{2mm}

(1) For a group $G$ and an element $s$ of $G$, we denote $\{ g \in G \, | \, sg = gs \}$ by $G^s$.
\vspace{1mm}

(2) For a transformation group $G$ of a space $M$, the isotropy subgroup of $G$ at $m_1, \cdots, m_k \in M$ is denoted by $G_{m_1, \cdots, m_k} = \{ g \in G \; | \; gm_1 = m_1, \cdots, gm_k = m_k \vspace{1mm}
\}$.

(3) For a $\R$-vector space $V$, its complexification $\{ u + iv \; | \; u, v \in V \}$ is denoted by $V^C$. The complex conjugation in $V^C$ is denoted by $\tau$ : $\tau(u + iv) = u - iv$. In particular, the complexification of $\R$ is briefly denoted by $C$ : $\R^C = C$.
\vspace{1mm}

(4) For a Lie group $G$, the Lie algebra of $G$ is denoted by the corresponding German small letter $\gg$. For example, $\so(n)$ is the Lie algebra of the group $SO(n)$.
\vspace{1mm}

(5) Although we will give all definitions used in the following Sections, if in case of insufficiency, refer to [3],[4] or [5]. 
\vspace{4mm}

\begin{center}
 {\large {\bf 1. Group $F_4$}}
\end{center}
\vspace{3mm}

We use the same notation as in [1], [2] or [5] (however, some will be rewritten). For example,
the Cayley algebra $\gC = \H \oplus \H e_4$,
\vspace{1mm}

the exceptional Jordan algebra $\gJ = \{X \in M(3, \gC) \, | \, X^* = X \}$, the Jordan multiplication $X \circ Y$, the inner product $(X, Y)$ and the elements $E_1, E_2, E_3 \in \gJ$,
\vspace{1mm}

the group $F_4 = \{\alpha \in \Iso_{\sR}(\gJ) \, | \, \alpha(X \circ Y) = \alpha X \circ \alpha Y \}$.\\

We define $\R$-linear transformations $\sigma$ and $\gamma$ of $\gJ$ by
$$
    \sigma X = \sigma\pmatrix{\xi_1 & x_3 & \ov{x}_2 \cr
                              \ov{x}_3 & \xi_2 & x_1 \cr 
                              x_2 & \ov{x}_1 & \xi_3} 
             = \pmatrix{\xi_1 &\!\! -x_3 &\!\! -\ov{x}_2 \cr
                        -\ov{x}_3 &\!\! \xi_2 &\!\! x_1 \cr 
                         -x_2 &\!\! \ov{x}_1 &\!\! \xi_3}, \;
     \gamma X = \pmatrix{\xi_1 &\!\! \gamma x_3 &\!\! \ov{\gamma x}_2 \cr
                              \ov{\gamma x}_3 &\!\! \xi_2 &\!\! {\gamma x}_1 \cr 
                              {\gamma x}_2 &\!\! \ov{\gamma x}_1 &\!\! \xi_3}, $$
respectively, where ${\gamma x}_k = \gamma(m_k + a_k e_4) = m_k - a_k e_4,\,\,x_k = m_k + a_k e_4 \in \H \oplus \H{e_4} = \gC.$ Then, $\sigma, \gamma \in F_4$ and $\sigma^2 = {\gamma}^2 = 1$. $\sigma$ and $\gamma$ are commutative. From $\sigma\gamma = \gamma\sigma$, we have
$$
              (F_4)^\sigma \cap (F_4)^{\gamma} = ((F_4)^{\sigma})^{\gamma} = ((F_4)^{\gamma})^{\sigma} .$$
Hence, this group will be denoted briefly by $(F_4)^{\sigma,\gamma}$.
\vspace{3mm}

{\sc Proposition 1.1.} \quad $(F_4)^\gamma \cong (Sp(1)\times Sp(3))/\Z_2, \,\, \Z_2 = \{(1, E), (-1, -E)\}$.
\vspace{2mm}

{\sc Proof.} The isomorphism is induced by the homomorphism $\varphi : Sp(1)\times Sp(3) \to (F_4)^\gamma,\,\,\varphi(p, A)(M + \a) = AMA^* + p\a{A^*},\,\,M + \a \in \gJ(3, \H) \oplus {\H}^3 = \gJ $. (In detail, see [3], [5].)
\vspace{3mm}

{\sc Lemma 1.2.} {\it $\varphi : Sp(1)\times Sp(3) \to (F_4)^\gamma $ of Proposition $1.1$ satisfies $\sigma\varphi(p, A)\sigma = \varphi(p, {I_1}A{I_1})$, where $I_1 = \diag(-1, 1, 1)$.}
\vspace{2mm}

{\sc Proof.} From $\sigma = \varphi(-1, I_1)$, we have the required one.
\vspace{3mm}

Now, we shall determine the group structure of $(F_4)^{\sigma,\gamma} = {((F_4)^\gamma)}^\sigma = {((F_4)^\sigma)}^\gamma = $ $(F_4)^\sigma \cap (F_4)^\gamma$.
\vspace{3mm}

{\sc Theorem 1.3.}  $(F_4)^{\sigma, \gamma} \cong (Sp(1)\times Sp(1)\times Sp(2))/\Z_2,\,\,\Z_2 = \{(1, 1, E), (-1, -1,$ $ -E)\}$.
\vspace{2mm}

{\sc Proof.} We define a map $\varphi_4 : Sp(1)\times Sp(1)\times Sp(2) \to (F_4)^{\sigma, \gamma}$ by 
$$
 {\varphi_4} (p, q, B )(M + \a) 
       = \left(\begin{array}{c|c}
            \!\!  q & \!\! \begin{array}{cc}
                     0 & \!\! 0
                   \end{array} \\
               \hline
              \!\! \begin{array}{c}
                 0 \\ 0
               \end{array} & \!\! B
               \end{array}  \right)  M \left(\begin{array}{c|c}
                                      \!\!  q & \!\! \begin{array}{cc}
                                                     0 & \!\!0
                                                  \end{array} \\
                                             \hline
                                             \!\! \begin{array}{c}
                                                0 \\ 0
                                             \end{array} & \!\! B
                                             \end{array}  \right)^{*} +  p\a\left(\begin{array}{c|c}
                                                                          \!\!   q & \!\! \begin{array}{cc}
                                                                                       0 & \!\! 0
                                                                                     \end{array} \\
                                                                                     \hline
                                                                                \!\! \begin{array}{c}
                                                                                     0 \\ 0
                                                                                 \end{array} & \!\! B
                                                                                 \end{array}  \right)^{*},
$$
$M + \a \in \gJ(3, \H)\oplus \H^3 = \gJ$, as the restriction of Proposition 1.1. By Lemma 1.2, $\varphi_4$ is well-defined and a homomorphism. We shall show that $\varphi_4$ is onto. Let $\alpha \in (F_4)^{\sigma, \gamma}$. Since $(F_4)^{\sigma, \gamma} \subset (F_4)^\gamma$, there exist $p \in Sp(1)$ and $A \in Sp(3)$ such that $\alpha = \varphi(p, A)$ (Proposition 1.1). From $\sigma \alpha \sigma = \alpha$, we have $\varphi(p, {I_1}A{I_1}) = \varphi(p, A)$ (Lemma 1.2). Hence, 
$$
\left\{\begin{array}{l}
       p = p
\vspace{1mm}\\
       {I_1}A{I_1} = A   \end{array} \right.
\quad \mbox{or} \quad
\left\{\begin{array}{l}
       p = -p
\vspace{1mm}\\
       {I_1}A{I_1} = -A   \end{array} \right. $$
The latter case is impossible because $p = 0$ is false. In the former case, from ${I_1}A{I_1} = A$,  we have 
\vspace{1mm}\\
$A = 
\left(\begin{array}{c|c}
            \!\!  q & \!\! \begin{array}{cc}
                     0 & \!\! 0
                   \end{array} \\
               \hline
              \!\! \begin{array}{c}
                 0 \\ 0
               \end{array} & \!\! B
               \end{array}  \right),\,\, q \in Sp(1), \,\,B \in Sp(2).$ Hence,
\vspace{1mm}\\
 $\alpha =$
$\varphi(\begin{array}{ll}
                                \!\! q, & \!\! \left(\begin{array}{c|c}
                                            \!\!  q & \!\! \begin{array}{cc}
                                                         0 & \!\! 0
                                                         \end{array} \\
                                                         \hline
                                                       \!\! \begin{array}{c}
                                                         0 \\ 0
                                                         \end{array} & \!\! B
                                                              \end{array}  \right)\end{array}) = {\varphi_4} (p, q, B)$,
\vspace{1mm}\\
 that is, $\varphi_4$ is onto. And $\Ker\varphi_4 = \{(1, 1, E), (-1,$
$-1, -E)\} = \Z_2$. Thus, we have the required isomorphism $(Sp(1)\times Sp(1)\times Sp(2))/\Z_2 \cong (F_4)^{\sigma, \gamma}$.
\vspace{3mm}

\begin{center}
{\large {\bf 2. Group $E_6$}}
\end{center}
\vspace{3mm}

We use the same notation as in [1], [2] or [5] (however, some will be rewritten). For example,
the complex exceptional Jordan algebra $\gJ^C = \{X \in M(3, \gC^C) \, | \, X^* = X \}$, the Freudenthal multiplication $X \times Y$ and the Hermitian inner product $\langle X, Y \rangle$,
\vspace{-3mm}

the group $E_6 = \{\alpha \in \Iso_C(\gJ^C) \, | \, \alpha X \times \alpha Y = \tau\alpha\tau(X \times Y), \langle \alpha X, \alpha Y \rangle = \langle X, Y \rangle \}$, and the natural inclusion $F_4 \subset E_6$.
\vspace{3mm}

{\sc Proposition 2.1.}\quad $(E_6)^\gamma \cong (Sp(1)\times SU(6))/\Z_2,\,\,\Z_2 = \{(1, E), (-1, -E)\}$.
\vspace{2mm}

{\sc Proof.} The isomorphism is induced by the homomorphism $\varphi : Sp(1) \times SU(6) \to (E_6)^\gamma,\,\,
  \varphi(p, A)(M + \a) = {k_J}^{-1}(A{k_J}(M){}^{t}A) + p\a{k^{-1}}(A^{*}), \, M + \a \in \gJ(3, \H)^C \oplus (\H^3)^C = \gJ^C $. (In detail, see [3], [5].)
\vspace{3mm}

{\sc Lemma 2.2.} {\it $\varphi : Sp(1)\times SU(6) \to (E_6)^\gamma$ of Proposition $2.1$ satisfies $\sigma\varphi(p, A)\sigma \\
= \varphi(p, {I_2}A{I_2})$, where $I_2 = \diag(-1, -1 ,1, 1, 1, 1)$.}
\vspace{2mm}

{\sc Proof.} From $\sigma = \varphi(-1, I_2)$, we have the required one.
\vspace{3mm}

Now, we shall determine the group structure of $(E_6)^{\sigma,\gamma} = {((E_6)^\gamma)}^\sigma = {((E_6)^\sigma)}^\gamma\\
 = (E_6)^\sigma \cap (E_6)^\gamma$.
\vspace{3mm}

{\sc Theorem 2.3.} $(E_6)^{\sigma,\gamma} \cong (Sp(1)\times S(U(2)\times U(4)))/\Z_2,\,\,\Z_2 = \{(1, E), (-1,$ $ -E)\}$.
\vspace{2mm}

{\sc Proof.} We define a map $\varphi_6 : Sp(1)\times S(U(2)\times U(4)) \to (E_6)^{\sigma,\gamma}$ by
$$
{ \varphi_6} (p, A)(M + \a) = {k_J}^{-1}(A{k_J}(M){}^{t}A) + p\a{k^{-1}}(A^{*}),$$
 $\, M + \a \in \gJ(3, \H)^C \oplus (\H^3)^C = \gJ^C $,
 as the restriction of $\varphi$ of Proposition 2.1. By Lemma 2.2, $\varphi_6$ is well-defined and a homomorphism. We shall show that $\varphi_6$ is onto. Let $\alpha \in (E_6)^{\sigma,\gamma}$. Since $(E_6)^{\sigma, \gamma} \subset (E_6)^\gamma$, there exist $p \in Sp(1)$ and $A \in SU(6)$ such that $\alpha = \varphi(p, A)$ (Proposition 2.1). From $\sigma \alpha \sigma = \alpha$, we have $\varphi(p, {I_2}A{I_2}) = \varphi(p, A)$ (Lemma 2.2). Hence, 
$$
\left\{\begin{array}{l}
       p = p
\vspace{1mm}\\
       {I_2}A{I_2} = A   \end{array} \right.
\quad \mbox{or} \quad
\left\{\begin{array}{l}
       p = -p
\vspace{1mm}\\
       {I_2}A{I_2} = -A   \end{array} \right. $$
The latter case is impossible because $p = 0$ is false. In the former case, we have $A \in S(U(2)\times U(4))$. Therefore,
$\varphi_6$ is onto. $\Ker\varphi_6 = \{(1, E), (-1, -E)\} = \Z_2$. Thus, we have
 the required isomorphism $(Sp(1)\times S(U(2)\times U(4)))/\Z_2 \cong (E_6)^{\sigma,\gamma}$.
\vspace{3mm}

\begin{center}
{\large {\bf 3. Group $E_7$}}
\end{center}
\vspace{3mm}

We use the same notation as in [1],[4] or [5] (however, some will be rewritten). For example,
the Freudenthal $C$-vector space $\gP^C = \gJ^C \oplus \gJ^C \oplus C \oplus C$, the Hermitian inner product $\langle P, Q \rangle$, the $C$-linear map $P \times Q$ : $\gP^C \to \gP^C$($P, Q \in \gP^C$),
\vspace{1mm}

the group $E_7 = \{\alpha \in \Iso_C(\gP^C) \, | \, \alpha(X \times Y)\alpha^{-1} = \alpha P \times \alpha Q, \langle \alpha P, \alpha Q \rangle = \langle P, Q \rangle \}$, the natural inclusion $E_6 \subset E_7$ and elements $\sigma, \sigma' \in F_4 \subset E_6 \subset E_7$, $\lambda \in E_7$.
\vspace{2mm}

We shall consider the following subgroup of $F_4$. 
$$
  ((F_4)^{\sigma,\gamma})_{F_1(h)} = \{\alpha \in (F_4)^{\sigma,\gamma}\, |\, \alpha F_1(h) = F_1(h) \mbox{ for all } h \in \H \} .$$
\vspace{-4mm}

{\sc Proposition 3.1.} \qquad $((F_4)^{\sigma,\gamma})_{F_1(h)} \cong Sp(1)\times Sp(1) (= Spin(4))$.
\vspace{2mm}

{\sc Proof.} We define a map $\varphi : Sp(1)\times Sp(1) \to ((F_4)^{\sigma,\gamma})_{F_1(h)}$ by
$$
\varphi (p, q)(M + \a) 
       = \left(\begin{array}{c|c}
            \!\!  q & \!\! \begin{array}{cc}
                     0 & \!\! 0
                   \end{array} \\
               \hline
              \!\! \begin{array}{c}
                 0 \\ 0
               \end{array} & \!\! E
               \end{array}  \right)  M \left(\begin{array}{c|c}
                                      \!\!  q & \!\! \begin{array}{cc}
                                                     0 & \!\!0
                                                  \end{array} \\
                                             \hline
                                             \!\! \begin{array}{c}
                                                0 \\ 0
                                             \end{array} & \!\! E
                                             \end{array}  \right)^{*} + p\a\left(\begin{array}{c|c}
                                                                          \!\!   q & \!\! \begin{array}{cc}
                                                                                       0 & \!\! 0
                                                                                     \end{array} \\
                                                                                     \hline
                                                                                \!\! \begin{array}{c}
                                                                                     0 \\ 0
                                                                                 \end{array} & \!\! E
                                                                                 \end{array}  \right)^{*},
$$
as the restriction of $\varphi_4$ of Theorem 1.3. By 
$F_1(h) = \pmatrix{0 & 0 & 0 \cr
                   0 & 0 & h \cr 
                   0 & \ov{h} & 0} + \O$, 
$\varphi$ is well-defined and homomorphism. We shall show that $\varphi$ is onto. Let $\alpha \in ((F_4)^{\sigma,\gamma})_{F_1(h)}$. Since $((F_4)^{\sigma,\gamma})_{F_1(h)} \subset (F_4)^{\sigma,\gamma}$, there exist $p, q \in Sp(1)$ and $B \in Sp(2)$ such that $\alpha = {\varphi_4}(p, q, B)$(Theorem 1.3). From $\alpha F_1(h) = F_1(h)$, we have
          $B\pmatrix{0 & h \cr 
                     \ov{h} & 0}B^* = \pmatrix{0 & h \cr 
                                             \ov{h} & 0}$, so that 
$$
   \alpha = {\varphi_4}(p, q, E) \quad \mbox{or} \quad \alpha = {\varphi_4}(p, q, -E).  $$
In the former case, we have $\alpha = {\varphi_4}(p, q, E) = \varphi(p, q)$. In the latter case, we have
\vspace{-5mm}

\begin{eqnarray*}
  \alpha \!\!\! &=& \!\!\!{\varphi_4}(p, q, -E) = {\varphi_4}(-p, -q, E) {\varphi_4}(-1, -1, -E)
\vspace{1mm}\\
         \!\!\! &=& \!\!\!{\varphi_4}(-p, -q, E){}1 = \varphi(-p, -q).
\end{eqnarray*}
\vspace{-5mm}

\noindent Hence, $\varphi$ is onto. $\Ker\varphi = \{(1, 1)\}$. Thus, we have the required isomorphism $Sp(1)\times Sp(1) \cong ((F_4)^{\sigma,\gamma})_{F_1(h)}$.
\vspace{2mm}

Hereafter, in $\gP^C$, we use the following notations. 
$$
\begin{array}{ll}  
    (F_1(h),0, 0, 0) = \dot{F_1}(h), \quad (0, E_1, 0, 1) = \ti{E}_1, 
\vspace{1mm}\\
    (0, E_1, 0, -1) = \ti{E}_{-1} , \quad (E_2 + E_3, 0, 0, 0) = \dot{E}_{23}.
\end{array}$$
\vspace{-2mm}

 We shall consider a subgroup $({((E_7)^{\kappa, \mu})}^\gamma)_{\dot{F_1}(h),\ti{E}_1,\ti{E}_{-1}, \dot{E}_{23}}$
of $E_7$. 
\vspace{3mm}

{\sc Lemma 3.2.} {\it The Lie algebra $({((\gge_7)^{\kappa, \mu})}^\gamma)_{\dot{F_1}(h),\ti{E}_1,\ti{E}_{-1},  \dot{E}_{23}}$ of the group \\

$({((E_7)^{\kappa, \mu})}^\gamma)_{\dot{F_1}(h),\ti{E}_1,\ti{E}_{-1},  \dot{E}_{23}}$ is given by} 
$$
\begin{array}{l}
({((\gge_7)^{\kappa, \mu})}^\gamma)_{\dot{F_1}(h),\ti{E}_1,\ti{E}_{-1},  \dot{E}_{23}}
\vspace{1.5mm}\\
= \Big \{ {\mit\Phi}(\left(\begin{array}{c|c}
                                0 & 0 \\
                             \hline
                                0 & D'_4
                            \end{array}   \right), 0, 0, 0)\, \Big |\, \left(\begin{array}{c|c}
                                                                           0 & 0 \\
                                                                       \hline
                                                                           0 & D'_4
                                                         \end{array}   \right) \in \so(8),\,\, D'_4 \in \so(4)\Big \}.
\end{array}$$
{\it In particular, we have} 
$$
  \dim(({((\gge_7)^{\kappa, \mu})}^\gamma)_{\dot{F_1}(h),\ti{E}_1,\ti{E}_{-1},  \dot{E}_{23}}) = 6.$$
\vspace{-3mm}

Hereafter,  $\left(\begin{array}{c|c}
                                0 & 0 \\
                             \hline
                                0 & D'_4
                            \end{array}   \right)$ will be denoted by $D'_4$, and also ${\mit\Phi}(D'_4, 0, 0, 0)$ will be denoted by ${\mit\Phi}_4$.
\vspace{2mm}

{\sc Proposition 3.3.}\quad $({((E_7)^{\kappa, \mu})}^\gamma)_{\dot{F_1}(h),\ti{E}_1,\ti{E}_{-1},  \dot{E}_{23}} = ((F_4)^{\sigma, \gamma})_{F_1(h)}$.
\vspace{2mm}

{\sc Proof.} Let $\alpha \in ((F_4)^{\sigma, \gamma})_{F_1(h)}$. Since $((F_4)^{\sigma, \gamma})_{F_1(h)} \subset (F_4)^\sigma = (F_4)_{E_1}$ (as for $(F_4)^\sigma = (F_4)_{E_1}$, see [3], [5]), we see $\alpha E_1 = E_1$. As a result, because $\kappa$ and $\mu$ are defined using by $E_1$ (see [1], [4] or [5]), we see that $\kappa\alpha = \alpha\kappa$ and $ \mu\alpha = \alpha\mu$. From $\alpha E = E$ \, (see [3], [5]), we have $\alpha(E_2 + E_3) = E_2 + E_3$. Hence, $\alpha\dot{E}_{23} = \dot{E}_{23}$. Moreover, from $\alpha (0, 0, 0, 1) = (0, 0, 0, 1)$ (see [4], [5]), we have $\alpha \ti{E}_1 = \ti{E}_1$ and $\alpha \ti{E}_{-1} = \ti{E}_{-1} $. Obviously $\alpha\dot{F_1}(h) = \dot{F_1}(h)$. Thus, $\alpha \in ({((E_7)^{\kappa, \mu})}^\gamma)_{\dot{F_1}(h),\ti{E}_1,\ti{E}_{-1},  \dot{E}_{23}}$. Conversely, let $\alpha \in 
({((E_7)^{\kappa, \mu})}^\gamma)_{\dot{F_1}(h),\ti{E}_1,\ti{E}_{-1},  \dot{E}_{23}}$. From $\alpha \ti{E}_1 = \ti{E}_1 $ and $\alpha\ti{E}_{-1} = \ti{E}_{-1}$, we have $\alpha (0, E_1, 0, 0) = (0, E_1, 0, 0)$ and $\alpha (0, 0, 0, 1) = (0, 0, 0, 1)$. Hence, $\alpha \in ((E_6)^\gamma)_{F_1(h), E_1, E_2 + E_3}$(see [4], [5]). Thus, $({((F_4)_{E_1})}^\gamma)_{F_1(h)} = ((F_4)^{\sigma, \gamma})_{F_1(h)}$. Therefore, the proof of this proposition is completed.
\vspace{2mm}

Next, we shall consider the following subgroup of $F_4$.
$$
  ((F_4)^{\sigma,\gamma})_{F_1(he_4)} = \{\alpha \in (F_4)^{\sigma,\gamma}\, |\, \alpha F_1(he_4) = F_1(he_4)  \mbox{ for all } h \in \H \} .$$
\vspace{-4mm}

{\sc Proposition 3.4.}\quad $((F_4)^{\sigma,\gamma})_{F_1(h{e_4})} \cong Sp(2) (= Spin(5))$.
\vspace{2mm}

{\sc Proof.} We define a map $\varphi : Sp(2) \to ((F_4)^{\sigma,\gamma})_{F_1(h{e_4})}$ by
$$
\varphi (B)(M + \a) 
       = \left(\begin{array}{c|c}
            \!\!  1 & \!\! \begin{array}{cc}
                     0 & \!\! 0
                   \end{array} \\
               \hline
              \!\! \begin{array}{c}
                 0 \\ 0
               \end{array} & \!\! B
               \end{array}  \right)  M \left(\begin{array}{c|c}
                                      \!\!  1 & \!\! \begin{array}{cc}
                                                     0 & \!\!0
                                                  \end{array} \\
                                             \hline
                                             \!\! \begin{array}{c}
                                                0 \\ 0
                                             \end{array} & \!\! B
                                             \end{array}  \right)^{*} + \a\left(\begin{array}{c|c}
                                                                          \!\!   1 & \!\! \begin{array}{cc}
                                                                                       0 & \!\! 0
                                                                                     \end{array} \\
                                                                                     \hline
                                                                                \!\! \begin{array}{c}
                                                                                     0 \\ 0
                                                                                 \end{array} & \!\! B
                                                                                 \end{array}  \right)^{*},
$$
as the restriction of $\varphi_4$ of Theorem 1.3. Obviously $\varphi$ is well-defined and homomorphism. We shall show that $\varphi$ is onto. Let $\alpha \in ((F_4)^{\sigma,\gamma})_{F_1(he_4)}$. Since $((F_4)^{\sigma,\gamma})_{F_1(h{e_4})} \subset (F_4)^{\sigma,\gamma}$, there exist $p, q \in Sp(1)$ and $B \in Sp(2)$ such that $\alpha = {\varphi_4}(p, q, B)$(Theorem 1.3). From $\alpha F_1(h{e_4}) = F_1(h{e_4}) (= O + (h, 0, 0))$, we have 
$ph\ov{q} = h(h \in \H)$, so that
$$
   \alpha = {\varphi_4}(1, 1, B) \quad \mbox{or} \quad \alpha = {\varphi_4}(-1, -1, B).  $$
In the former case, we have $\alpha = {\varphi_4}(1, 1, B) = \varphi(B)$. In the latter case, we have
\vspace{-5mm}

\begin{eqnarray*}
  \alpha \!\!\! &=& \!\!\!{\varphi_4}(-1, -1, B) = {\varphi_4}(1, 1, -B) {\varphi_4}(-1, -1, -E)
\vspace{1mm}\\
         \!\!\! &=& \!\!\!{\varphi_4}(1, 1, -B){}1 = \varphi(-B).
\end{eqnarray*}
\vspace{-5mm}

\noindent Hence, $\varphi$ is onto. $\Ker\varphi = \{E \}$. Thus, we have the required isomorphism $Sp(2) \cong ((F_4)^{\sigma,\gamma})_{F_1(h{e_4})}$.
\vspace{2mm}

Then, we have the following proposition. 
\vspace{2mm}

{\sc Proposition 3.5.}\quad $({((E_7)^{\kappa, \mu})}^\gamma)_{\dot{F_1}(h{e_4}),\ti{E}_1,\ti{E}_{-1}, \dot{E}_{23}} = ((F_4)^{\sigma, \gamma})_{F_1(h{e_4})}$.
\vspace{2mm}

{\sc Proof.} This proof is in the way similar to Proposition 3.3.
\vspace{3mm}

We shall consider the subgroup
$({((E_7)^{\kappa, \mu})}^\gamma)_{\dot{F_1}(h{e_4}),\ti{E}_1,\ti{E}_{-1}}$ of $E_7$.
\vspace{2mm}

{\sc Lemma 3.6.}  {\it The Lie algebra $({((\gge_7)^{\kappa, \mu})}^\gamma)_{\dot{F_1}(h{e_4}),\ti{E}_1,\ti{E}_{-1}}$ of the group \\
 $({((E_7)^{\kappa, \mu})}^\gamma)_{\dot{F_1}(h{e_4}),\ti{E}_1,\ti{E}_{-1}}$ is given by} 
$$
\begin{array}{l}
({((\gge_7)^{\kappa, \mu})}^\gamma)_{\dot{F_1}(h{e_4}),\ti{E}_1,\ti{E}_{-1}}
\vspace{1mm}\\
= \Big \{ {\mit\Phi} \Bigl( \left(\begin{array}{c|c}
                                D_4 & 0 \\
                             \hline
                                0 & 0
                            \end{array} \right) + \wti{A}_1(p) + i\pmatrix{0 & 0 & 0 \cr
                                                                               0 & \epsilon & q \cr
                                                                               0 & \ov{q} & -\epsilon}^{\sim}, 
                                 0, 0, 0 \Bigr)\, \Big |\, \left( \begin{array}{c|c}
                                                                           D_4 & 0 \\
                                                                       \hline
                                                                           0 & 0
                                                         \end{array}   \right) \in \so(8),
\vspace{2mm}\\
 \qquad \qquad \qquad \qquad \qquad \qquad \qquad                              D_4 \in \so(4),\, \epsilon \in \R,\, p, q \in \H \Big \}.
\end{array}$$
{\it In particular, we have} 
$$
  \dim(({((\gge_7)^{\kappa, \mu})}^\gamma)_{\dot{F_1}(h{e_4}),\ti{E}_1,\ti{E}_{-1}}) = 15.$$
 
Hereafter, $\left(\begin{array}{c|c}
                                D_4 & 0 \\
                             \hline
                                0 & 0
                            \end{array}   \right)$ will be denoted by $D_4$.
\vspace{3mm}

{\sc Lemma 3.7.}\quad (1) {\it  For $a \in \H$, we define a map ${\wti{\alpha}}_1(a)$ of $\gJ^C$ by}
\vspace{-5mm}\\
 
$$
\left \{ \begin{array}{l}
      \xi'_1 = \xi_1 
\vspace{1mm}\\
      \xi'_2 = \dfrac{\xi_2 - \xi_3}{2} + \dfrac{\xi_2 + \xi_3}{2}\cos|a| + i\dfrac{(a, x_1)}{|a|}\sin|a| 
\vspace{1mm}\\
      \xi'_3 = -\dfrac{\xi_2 - \xi_3}{2} + \dfrac{\xi_2 + \xi_3}{2}\cos|a| + i\dfrac{(a, x_1)}{|a|}\sin|a| 
\end{array} \right . $$
$$
\,\,\,\,\left \{ \begin{array}{l}
      x'_1 = x_1 + i\dfrac{(\xi_2 + \xi_3)a}{|a|}\sin|a| - {\dfrac{2(a, x_1)a}{|a|^2}}{(\sin{\dfrac{|a|}{2}})^2}
\vspace{1mm}\\
      x'_2 = x_2\cos \dfrac{|a|}{2} + i\dfrac{\ov{x_3 a}}{|a|} \sin \dfrac{|a|}{2}
\vspace{1mm}\\
      x'_3 = x_3\cos \dfrac{|a|}{2} + i\dfrac{\ov{a x_2}}{|a|} \sin \dfrac{|a|}{2}
\end{array} \right. $$
{\it Then}, $\wti{\alpha}_1(a) \in ({((E_7)^{\kappa, \mu})}^\gamma)_{\dot{F_1}(h{e_4}),\ti{E}_1,\ti{E}_{-1}}$. \\
\vspace{-2mm}

(2) {\it For $t \in \R$, we define a map $\wti{\alpha}_{23}(t)$ of $\gJ^C$ by}

$$
\wti{\alpha}_{23}(t)\pmatrix{\xi_1 & x_3 & \ov{x}_2 \cr
                  \ov{x}_3 & \xi_2 & x_1 \cr
                  x_2 & \ov{x}_1 & \xi_3} 
                = \pmatrix{\xi_1 & e^{it/2}x_3 & e^{-it/2}\ov{x}_2 \cr
                           e^{it/2}\ov{x}_3 & e^{it}\xi_2 & x_1 \cr
                           e^{-it/2}x_2 & \ov{x}_1 & e^{-it}\xi_3}. $$
{\it Then}, $\wti{\alpha}_{23}(t) \in  ({((E_7)^{\kappa, \mu})}^\gamma)_{\dot{F_1}(h{e_4}),\ti{E}_1,\ti{E}_{-1}}$.
\vspace{3mm}

{\sc Proof.}(1) For $a \in \H$, we have $i\wti{F}_1 (a) \in ({((\gge_7)^{\kappa, \mu})}^\gamma)_{\dot{F_1}(h{e_4}),\ti{E}_1,\ti{E}_{-1}}$(Lemma 3.6). Hence, $\wti{\alpha}_1 (a) = \exp{i\wti{F}_1 (a)} \in ({((E_7)^{\kappa, \mu})}^\gamma)_{\dot{F_1}(h{e_4}),\ti{E}_1,\ti{E}_{-1}}$. \\ 
\vspace{-2mm}

(2) For $t \in \R$, we have $it(E_2-E_3)^{\sim} \in ({((\gge_7)^{\kappa, \mu})}^\gamma)_{\dot{F_1}(h{e_4}),\ti{E}_1, \ti{E}_{-1}}$(Lemma 3.6). Hence, $\wti{\alpha}_{23} (t) = \exp{it(E_2-E_3)^{\sim}} \in ({((E_7)^{\kappa, \mu})}^\gamma)_{\dot{F_1}(h{e_4}),\ti{E}_1,\ti{E}_{-1}}$.
\vspace{3mm}

We define a $6$ dimensional $\R$-vector space $V^6$ by
\begin{eqnarray*}
    V^6 \!\!\!& = &\!\!\! \{P \in \gP^C \,|\,\kappa P = P, \mu\tau\lambda P = P, \gamma P = P, \langle P, \ti{E_1} \rangle = 0, \langle P, \ti{E}_{-1} \rangle = 0 \} 
\vspace{1mm}\\
        \!\!\!& = &\!\!\! \Big \{P = \Bigl(\pmatrix{0 & 0 & 0 \cr
                                              0 & \xi & h \cr
                                              0 & \ov{h} & -\tau{\xi}}, 0, 0, 0 \Bigr)\,\Big |\,\xi \in C,\,\,h \in \H \Big \}
\end{eqnarray*}
with the norm (see [5] for the definition of $\{{\quad },{\quad }\}$'s)
$$ 
      (P, P)_{\mu} = \dfrac{1}{2}\{\mu P, P \} = \dfrac{1}{2}(\mu P, \lambda P) = (\tau\xi)\xi + \ov{h}h. $$

\noindent Then, $S^5 = \{P \in V^6\,|\, (P, P)_{\mu} = 1 \}$ is a $5$ dimensional sphere.
\vspace{3mm}

{\sc Lemma 3.8}\qquad $({((E_7)^{\kappa, \mu})}^\gamma)_{\dot{F_1}(h{e_4}),\ti{E}_1,\ti{E}_{-1}}/Spin(5) \simeq S^5$.\\
{\it In particular, $({((E_7)^{\kappa, \mu})}^\gamma)_{\dot{F_1}(h{e_4}),\ti{E}_1,\ti{E}_{-1}}$ is connected.}
\vspace{2mm}

{\sc Proof.} Since $E_7$ is commutative with $\tau \lambda$, the group $({((E_7)^{\kappa, \mu})}^\gamma)_{\dot{F_1}(h{e_4}),\ti{E}_1,\ti{E}_{-1}}$ acts on $S^5$. We shall show that this action is transitive. To show this, it is sufficient to show that any element $P \in S^5$ can be transformed to $(i(E_2+E_3), 0, 0, 0) \in S^5$ under the action of $({((E_7)^{\kappa, \mu})}^\gamma)_{\dot{F_1}(h{e_4}),\ti{E}_1,\ti{E}_{-1}}$. Now, for a given
$$
       P = \Bigl(\pmatrix{0 & 0 & 0 \cr
                     0 & \xi & h \cr
                     0 & \ov{h} & -\tau{\xi}}, 0, 0, 0 \Bigr) \in S^5,  $$
choose $t \in \R$ such that $e^{it}\xi \in \R$. For this $t \in \R$, operate $\wti{\alpha}_{23}(t)\,\mbox{(Lemma 3.7(2))} \in ({((E_7)^{\kappa, \mu})}^\gamma)_{\dot{F_1}(h{e_4}),\ti{E}_1,\ti{E}_{-1}}$ on $P$. Then, we have 
$$ 
     \wti{\alpha}_{23} (t)P =  \Bigl(\pmatrix{0 & 0 & 0 \cr
                                              0 & r & h \cr
                                              0 & \ov{h} & -r}, 0, 0, 0 \Bigr) = P_1,\,\, r \in \R.   $$
In the case of $h \ne 0$, operate $\wti{\alpha}_1 (\pi h/2|h|)\,\mbox{(Lemma 3.7(1))} \in ({((E_7)^{\kappa, \mu})}^\gamma)_{\dot{F_1}(h{e_4}),\ti{E}_1,\ti{E}_{-1}}$ on $P_1$. Then, we have
$$
    \wti{\alpha}_1 \bigl(\dfrac{\pi h}{2|h|} \bigr) P_1 = \Bigl(\pmatrix{0 & 0 & 0 \cr
                                              0 & \xi' & 0 \cr
                                              0 & 0 & -\tau\xi'}, 0, 0, 0 \Bigr) = P_2 \in S^5,\,\, \xi' \in C. $$
Here, from $(\tau\xi')\xi' = 1, \, \xi' \in C$, we can put $\xi' = e^{i\theta}, 0 \leq \theta < 2\pi$. Operate $\wti{\alpha}_{23} (-\theta)$ on $P_2$. Then, 
$$
     \wti{\alpha}_{23} \bigl(-\theta \bigr) P_2 = (E_2-E_3, 0, 0, 0) = P_3.  $$
Moreover, operate $\wti{\alpha}_{23} (\pi/2)$ on $P_3$,
$$
   \wti{\alpha}_{23} \bigl(\dfrac{\pi}{2} \bigr) P_3 = (i(E_2+E_3), 0, 0, 0) = i\dot{E}_{23} .  $$
This shows the transitivity. The isotropy subgroup $(((E_7)^{\kappa, \mu})^\gamma)_{\dot{F_1}(h{e_4}),\ti{E}_1,\ti{E}_{-1}}$ at $\dot{E}_{23}$ is $({((E_7)^{\kappa, \mu})}^\gamma)_{\dot{F_1}(h{e_4}),\ti{E}_1,\ti{E}_{-1},\dot{E}_{23}} = Sp(2)$ (Propositions 3.4, 3.5) $= Spin(5)$.
 Therefore, we have the homeomorphism \\
$(((E_7)^{\kappa, \mu})^\gamma)_{\dot{F_1}(h{e_4}),\ti{E}_1,\ti{E}_{-1}}/Spin(5) $ $ \simeq S^5$.
\vspace{3mm} 

{\sc Proposition 3.9.} \qquad $({((E_7)^{\kappa, \mu})}^\gamma)_{\dot{F_1}(h{e_4}),\ti{E}_1,\ti{E}_{-1}} \cong Spin(6)$.
\vspace{2mm}

{\sc Proof.} Since $({((E_7)^{\kappa, \mu})}^\gamma)_{\dot{F_1}(h{e_4}),\ti{E}_1,\ti{E}_{-1}}$ is connected (Lemma 3.8), we can define a homormorphism $
     \pi : ({((E_7)^{\kappa, \mu})}^\gamma)_{\dot{F_1}(h{e_4}),\ti{E}_1,\ti{E}_{-1}} \to SO(6) = SO(V^6)  $ by  
$$
       \pi(\alpha) = \alpha|V^6 .$$
It is not difficult to see that $\Ker\varphi = \{1, \sigma \} = \Z_2$. Since \\
 $\dim(({((E_7)^{\kappa, \mu})}^\gamma)_{\dot{F_1}(h{e_4}),\ti{E}_1,\ti{E}_{-1}}) = 15\,{\mbox{(Lemma 3.6)}}= \dim(\so(6))$, $\pi$ is onto.
 Hence, $({((E_7)^{\kappa, \mu})}^\gamma)_{\dot{F_1}(h{e_4}),\ti{E}_1,\ti{E}_{-1}}/\Z_2 \cong SO(6)$. Therefore, 
$({((E_7)^{\kappa, \mu})}^\gamma)_{\dot{F_1}(h{e_4}),\ti{E}_1,\ti{E}_{-1}}$ is isomorphism to $Spin(6)$ as a double covering group of $SO(6)$.
\vspace{2mm}

We shall consider a subgroup 
$({((E_7)^{\kappa, \mu})}^\gamma)_{\dot{F_1(h{e_4})},\ti{E}_1}$ of $E_7$. 
\vspace{2mm}

{\sc Lemma 3.10.} {\it The Lie algebra $({((\gge_7)^{\kappa, \mu})}^\gamma)_{\dot{F_1}(h{e_4}),\ti{E}_1}$ of the group \\
$({((E_7)^{\kappa, \mu})}^\gamma)_{\dot{F_1}(h{e_4}),\ti{E}_1}$ is given by }
$$
\begin{array}{l}
({((\gge_7)^{\kappa, \mu})}^\gamma)_{\dot{F_1}(h{e_4}),\ti{E_1}}
\vspace{1mm}\\
= \Big \{ {\mit\Phi}\Bigl( D_4  + \wti{A}_1(p) + i\pmatrix{0 & 0 & 0 \cr
                                                        0 & \epsilon & q \cr
                                                        0 & \ov{q} & -\epsilon}^{\sim}                                                                                                 , \pmatrix{0 & 0 & 0 \cr
                                                                           0 & \alpha & ix \cr
                                                                           0 & \ov{ix} & \tau\alpha}
, -\tau{\pmatrix{0 & 0 & 0 \cr                         
                 0 & \alpha & ix \cr                                                                            
                 0 & \ov{ix} & \tau\alpha}}, 0 \Bigr)
\end{array} $$
$$
\begin{array}{l} 
\vspace{1mm}\\
 \qquad \qquad \quad           \,\Big|\, D_4 \in \so(4) \subset \so(8), \epsilon \in \R,\alpha \in C, p, q, x \in \H \Big \}.
\end{array}  $$  {\it  In particular, we have }
$$
  \dim(({((\gge_7)^{\kappa, \mu})}^\gamma)_{\dot{F_1}(h{e_4}),\ti{E}_1}) = 21.$$
\vspace{-5mm}

{\sc Lemma 3.11.} {\it For $a \in \R$, we define maps $\alpha_k(a), k = 2, 3$ of $\gP^C$ by} 
$$
   \alpha_k(a)\pmatrix{X \vspace{1mm}\cr
                      Y \vspace{1mm}\cr
                      \xi \vspace{1mm}\cr
                      \eta}
   = \pmatrix{(1 + (\cos a - 1)p_k)X - 2(\sin a)E_k \times Y + \eta(\sin a)E_k \vspace{1mm}\cr
   2(\sin a)E_k \times X + (1 + (\cos a - 1)p_k)Y - \xi(\sin a)E_k 
\vspace{1mm}\cr
   ((\sin a)E_k, Y) + (\cos a)\xi
\vspace{1mm}\cr
   (- (\sin a)E_k, X) + (\cos a)\eta}, $$
{\it where $p_k : \gJ^C \to \gJ^C$ is defined by}
$$
   p_k(X) = (X, E_k)E_k + 4E_k \times (E_k \times X), \quad X \in \gJ^C.$$
{\it Then, $\alpha_k \in E_7$ and}
$ \alpha_2(a), \alpha_3(b) (a, b \in \R)$ {\it commute with each other.}
\vspace{2mm}

{\sc Proof.} For ${\mit\Phi}_k(a) = {\mit\Phi}(0, aE_k, -aE_k, 0) \in \gge_7$, we have $\alpha_k(a) = \exp{\mit\Phi}_k(a) \in E_7$. Since $[{\mit\Phi}_2(a), {\mit\Phi}_3(b)] = 0, \alpha_2(a)$ and $\alpha_3(b)$ are commutative.
\vspace{3mm}

We define a $7$ dimensional $\R$-vector space $V^7$ by
\begin{eqnarray*}
    V^7 \!\!\! & = & \!\!\! \{P \in \gP^C \,|\,\kappa P = P, \mu\tau\lambda P = P, \gamma P = P, \langle P, \ti{E_1}\rangle = 0 \} 
\vspace{1mm}\\
        \!\!\! & = & \!\!\! \Big \{P = \Bigl(\pmatrix{0 & 0 & 0 \cr
                                              0 & \xi & h \cr
                                              0 & \ov{h} & -\tau{\xi}}, \pmatrix{i\eta & 0 & 0 \cr
                                                                               0 & 0 & 0 \cr
                                                                               0 & 0 & 0},0, -i\eta \Bigr) 
   \,\Big |\,\xi \in C,\,\,h \in \H,\,\,\eta \in \R \Big \}
\end{eqnarray*}
with the norm 
$$ 
      (P, P)_{\mu} = \dfrac{1}{2}(\mu P, \lambda P) = (\tau\xi)\xi + \ov{h}h + \eta^2. $$
Then, $S^6 = \{P \in V^7\,|\, (P, P)_{\mu} = 1 \}$ is a $6$ dimensional sphere.
\vspace{3mm} 

{\sc Lemma 3.12}\qquad $({((E_7)^{\kappa, \mu})}^\gamma)_{\dot{F_1}(h{e_4}),\ti{E}_1}/Spin(6) \simeq S^6$.\\
{\it In particular, $({((E_7)^{\kappa, \mu})}^\gamma)_{\dot{F_1}(h{e_4}),\ti{E}_1}$ is connected.}
\vspace{2mm}

{\sc Proof.} The group $({((E_7)^{\kappa, \mu})}^\gamma)_{\dot{F_1}(h{e_4}),\ti{E}_1}$ acts on $S^6$. We shall show that this action is transitive. To show this, it is sufficient to show that any element $P \in S^6$ can be transformed to $(0, -iE_1, 0, i) \in S^6$ under the action of $({((E_7)^{\kappa, \mu})}^\gamma)_{\dot{F_1}(h{e_4}),\ti{E}_1}$. 
Now, for a given
$$
       P = \Bigl(\pmatrix{0 & 0 & 0 \cr
                     0 & \xi & h \cr
                     0 & \ov{h} & -\tau{\xi}}, \pmatrix{i\eta & 0 & 0 \cr
                                                          0 & 0 & 0 \cr
                                                          0 & 0 & 0}, 0, -i\eta \Bigr) \in S^6,  $$
choose $a \in \R, 0\leq a < \pi/2$ such that $\tan 2a = \dfrac{i2\eta}{\tau\xi-\xi}$(if $\tau\xi-\xi=0$, then 
let $a=\pi/4$). Operate $\alpha_{23}(a):= {\alpha_2}(a){\alpha_3}(a) = \exp({\mit\Phi}(0, a(E_2+E_3), -a(E_2+E_3), 0))\,\mbox{(Lemma 3.11)} \in ({((E_7)^{\kappa, \mu})}^\gamma)_{\dot{F_1}(h{e_4}),\ti{E}_1}$(Lemma 3.10) on $P$.
Then, the $\eta-$ \\
term of $\alpha_{23}(a)P$ is $(1/2)(\xi-\tau\xi)\sin 2a + i\eta\cos 2a = 0$. Hence, 
$$
   \alpha_{23} (a) P = \Bigl( \pmatrix{0 & 0 & 0 \cr
                            0 & \zeta & m \cr
                            0 & \ov{m} & -\tau\zeta}, 0, 0, 0 \Bigr) = P_1 \in S^5 \subset S^6. $$

\noindent Since $({((E_7)^{\kappa, \mu})}^\gamma)_{\dot{F_1}(h{e_4}),\ti{E}_1, \ti{E}_{-1}}( \subset ({((E_7)^{\kappa, \mu})}^\gamma)_{\dot{F_1}(h{e_4}),\ti{E_1}})$ acts transitivity on $S^5$ \\
(Lemma 3.8), there exist $\beta \in ({((E_7)^{\kappa, \mu})}^\gamma)_{\dot{F_1}(h{e_4}),\ti{E}_1, \ti{E}_{-1}}$ such that $\beta P_1 = (i(E_2+E_3), 0, 0, 0) = P_2 \in S^5 \subset S^6$. Moreover, operate $\alpha_{23} (-\pi/4)$
 on $P_2$,
$$
        \alpha_{23} \bigl(-\dfrac{\pi}{4} \bigr) P_2 = (0, -iE_1, 0, i) = -i\ti{E}_{-1}.   $$
This shows the transitivity. The isotropy subgroup $({((E_7)^{\kappa, \mu})}^\gamma)_{\dot{F_1}(h{e_4}),\ti{E}_1}$ at $\ti{E}_{-1}$ is $({((E_7)^{\kappa, \mu})}^\gamma)_{\dot{F_1}(h{e_4}),\ti{E}_1, \ti{E}_{-1}} = Spin(6)$(Proposition 3.9). Thus, we have the homeomorphism $({((E_7)^{\kappa, \mu})}^\gamma)_{\dot{F_1}(h{e_4}),\ti{E}_1}/Spin(6) \simeq S^6.$
\vspace{3mm}

{\sc Proposition 3.13.}\qquad $({((E_7)^{\kappa, \mu})}^\gamma)_{\dot{F_1}(h{e_4}),\ti{E}_1} \cong Spin(7).$
\vspace{2mm}

{\sc Proof.} Since $({((E_7)^{\kappa, \mu})}^\gamma)_{\dot{F_1}(h{e_4}),\ti{E}_1}$ is connected (Lemma 3.12), we can define a homormorphism $
     \pi : ({((E_7)^{\kappa, \mu})}^\gamma)_{\dot{F_1}(h{e_4}),\ti{E}_1} \to SO(7) = SO(V^7)  $ by  
$$
       \pi(\alpha) = \alpha|V^7 .$$
It is not difficult to see that $\Ker\varphi = \{1, \sigma \} = \Z_2$. Since \\
 $\dim(({((E_7)^{\kappa, \mu})}^\gamma)_{\dot{F_1}(h{e_4}),\ti{E}_1}) = 21\,{\mbox{(Lemma 3.10)}}= \dim(\so(7))$, $\pi$ is onto. Hence, \\
$({((E_7)^{\kappa, \mu})}^\gamma)_{\dot{F_1}(h{e_4}),\ti{E}_1}/\Z_2 \cong SO(7)$. Therefore, $({((E_7)^{\kappa, \mu})}^\gamma)_{\dot{F_1}(h{e_4}),\ti{E}_1}$ is isomorphism to $Spin(7)$ as a double covering group of $SO(7)$.
\vspace{2mm}

We shall consider the subgroup 
$({((E_7)^{\kappa, \mu})}^\gamma)_{\dot{F_1}(h{e_4})}$ of $E_7$.
\vspace{2mm}

{\sc Lemma 3.14.} {\it The Lie algebra $({((\gge_7)^{\kappa, \mu})}^\gamma)_{\dot{F_1}(h{e_4})}$ of the group \\
$({((E_7)^{\kappa, \mu})}^\gamma)_{\dot{F_1}(h{e_4})}$ is given by} 
$$
\begin{array}{l}
({((\gge_7)^{\kappa, \mu})}^\gamma)_{\dot{F_1}(h{e_4})}
\vspace{1mm}\\
= \Big \{ {\mit\Phi}\Bigl( D_4  + \wti{A}_1(p) + i\pmatrix{\epsilon_1 & 0 & 0 \cr
                                                               0 & \epsilon_2 & q \cr
                                                               0 & \ov{q} & \epsilon_3}^{\sim}                                                                                                 , \pmatrix{0 & 0 & 0 \cr
                                                                           0 & \alpha_2 & x \cr
                                                                           0 & \ov{x} & \alpha_3}
,         -\tau \pmatrix{0 & 0 & 0 \cr                         
                 0 & \alpha_2 & x \cr                                                                                            0 & \ov{x} & \alpha_3 },
\vspace{1mm}\\
 -\dfrac{3}{2}i{\epsilon_1} \Bigr) \,\Big|\, D_4 \in \so(4) \subset \so(8), \alpha_k \in C, p, q \in \H, x \in \H^C,\epsilon_k \in \R, \epsilon_1+\epsilon_2+\epsilon_3 = 0 \Big \}. 
\end{array}$$
{\it In particular, we have }
$$
  \dim(({((\gge_7)^{\kappa, \mu})}^\gamma)_{\dot{F_1}(h{e_4})}) = 28. $$
Hereafter, any element of the Lie algebra $({((\gge_7)^{\kappa, \mu})}^\gamma)_{\dot{F_1}(h{e_4})}$ will be  denoted by ${\mit\Phi}_8$.
\vspace{-1mm}

{\sc Lemma 3.15.} {\it For $t \in \R$, we define a map $\alpha(t)$ of $\gP^C$ by} 
$$
\begin{array}{l}
     \alpha(t)(X, Y, \xi, \eta) 
\vspace{1mm}\\
\quad
   = \Big(\pmatrix{e^{2it}\xi_1 & e^{it}x_3 & e^{it}\ov{x}_2 \cr
                    e^{it}\ov{x}_3 & \xi_2 & x_1 \cr
                    e^{it}x_2 & \ov{x}_1 & \xi_3}, 
           \pmatrix{e^{-2it}\eta_1 & e^{-it}y_3 & e^{-it}\ov{y}_2 \cr
                    e^{-it}\ov{y}_3 & \eta_2 & y_1 \cr
                    e^{-it}y_2 & \ov{y}_1 & \eta_3}, e^{-2it}\xi, e^{2it}\eta \Big).
\end{array}$$
{\it Then}, $\alpha(t) \in ({((E_7)^{\kappa, \mu})}^\gamma)_{\dot{F_1}(h{e_4})}.$
\vspace{2mm}

{\sc  Proof.} For ${\mit\Phi} = {\mit\Phi}(2itE_1 \vee E_1, 0, 0, -2it) \in ({((\gge_7)^{\kappa, \mu})}^\gamma)_{\dot{F_1}(h{e_4})}$ (Lemma 3.14), we have $\alpha(t) = \exp{\mit\Phi} \in ({((E_7)^{\kappa, \mu})}^\gamma)_{\dot{F_1}(h{e_4})}$ by $E_1 \vee E_1 = (1/3)(2E_1 - E_2 - E_3)^{\sim}$. 
\vspace{3mm}

 We define an $8$ dimensional $\R$-vector space $V^8$ by
\begin{eqnarray*}
    V^8 \!\!\! & = & \!\!\! \{P \in \gP^C \,|\,\kappa P = P, \mu\tau\lambda P = P, \gamma P = P \} \\
        \!\!\! & = & \!\!\! \Big \{P = \Bigl(\pmatrix{0 & 0 & 0 \cr
                                              0 & \xi & h \cr
                                              0 & \ov{h} & -\tau{\xi}}, \pmatrix{\eta & 0 & 0 \cr
                                                                               0 & 0 & 0 \cr
                                                                               0 & 0 & 0},0, \tau\eta \Bigr) 
                \,\Big |\,\xi,\eta \in C,\,\,h \in \H \Big \}
\end{eqnarray*}
with the norm 
$$ 
      (P, P)_{\mu} = \dfrac{1}{2}(\mu P, \lambda P) = (\tau\xi)\xi + \ov{h}h + (\tau\eta)\eta. $$
Then, $S^7 = \{ P \in V^8\,|\, (P, P)_{\mu} = 1 \}$ is a $7$ dimensional sphere.
\vspace{3mm}

{\sc Lemma 3.16.}\qquad $({((E_7)^{\kappa, \mu})}^\gamma)_{\dot{F_1}(h{e_4})}/Spin(7) \simeq S^7$.\\
{\it In particular, $({((E_7)^{\kappa, \mu})}^\gamma)_{\dot{F_1}(h{e_4})}$ is connected.}
\vspace{2mm}

{\sc Proof.} The group $({((E_7)^{\kappa, \mu})}^\gamma)_{\dot{F_1}(h{e_4})}$ acts on $S^7$. We shall show that this action is transitive. To show this, it is sufficient to show that any element $P \in S^7$ can be transformed to $(0, E_1, 0, 1) \in S^7$ under the action of $({((E_7)^{\kappa, \mu})}^\gamma)_{\dot{F_1}(h{e_4})}$. Now, for a given \\
$$
       P = \Bigl(\pmatrix{0 & 0 & 0 \cr
                     0 & \xi & h \cr
                     0 & \ov{h} & -\tau{\xi}}, \pmatrix{\eta & 0 & 0 \cr
                                                          0 & 0 & 0 \cr
                                                          0 & 0 & 0}, 0, \tau\eta \Bigr) \in S^7,  $$
choose $t \in \R$ such that $e^{-2it}\eta \in i\R$. Operate $\alpha(t)$ (Lemma 3.15)$\, \in ({((E_7)^{\kappa, \mu})}^\gamma)_{\dot{F_1}(h{e_4})}$ on $P$. Then, 
$$
       \alpha(t) P = \Bigl(\pmatrix{0 & 0 & 0 \cr
                     0 & \xi & h \cr
                     0 & \ov{h} & -\tau{\xi}}, \pmatrix{i\eta' & 0 & 0 \cr
                                                          0 & 0 & 0 \cr
                                                          0 & 0 & 0}, 0, -i\eta' \Bigr) = P_1 \in S^6 \subset S^7,\,\eta' \in \R   $$
Since $({((E_7)^{\kappa, \mu})}^\gamma)_{\dot{F_1}(h{e_4}), \ti{E}_1} (\subset ({((E_7)^{\kappa, \mu})}^\gamma)_{\dot{F_1}(h{e_4})})$ acts transitivity on $S^6$ \,(Lemma 3.12), there exists $\beta \in ({((E_7)^{\kappa, \mu})}^\gamma)_{\dot{F_1}(h{e_4}), \ti{E}_1}$ such that $\beta P_1 = (0, -iE_1, 0, i) = P_2 \in S^6 \subset S^7$. Moreover, operate $\alpha(-\pi/4)$ (Lemma 3.15) on $P_2$,
$$
       \alpha\bigl(-\dfrac{\pi}{4} \bigr) P_2 = (0, E_1, 0, 1) = \ti{E}_1.  $$
This shows the transitivity. The isotropy subgroup $({((E_7)^{\kappa, \mu})}^\gamma)_{\dot{F_1}(h{e_4})}$ at $\ti{E}_1$ is $({((E_7)^{\kappa, \mu})}^\gamma)_{\dot{F_1}(h{e_4}), \ti{E}_1} = Spin(7)$\,(Proposition 3.12). Thus, we have the homeomorphism $({((E_7)^{\kappa, \mu})}^\gamma)_{\dot{F_1}(h{e_4})}/Spin(7) \simeq S^7$.
\vspace{3mm}

{\sc Proposition 3.17.}\qquad $({((E_7)^{\kappa, \mu})}^\gamma)_{\dot{F_1}(h{e_4})} \cong Spin(8)$.
\vspace{2mm}

{\sc Proof.} Since $({((E_7)^{\kappa, \mu})}^\gamma)_{\dot{F_1}(h{e_4})}$ is connected (Lemma 3.16), we can define a homormorphism $
     \pi : ({((E_7)^{\kappa, \mu})}^\gamma)_{\dot{F_1}(h{e_4})} \to SO(8) = SO(V^8)  $ by  
$$
       \pi(\alpha) = \alpha|V^8 .$$
It is not difficult to see that $\Ker\varphi = \{1, \sigma \} = \Z_2$. Since $\dim(({((E_7)^{\kappa, \mu})}^\gamma)_{\dot{F_1}(h{e_4})}) = 28\,{\mbox{(Lemma 3.14)}}= \dim(\so(8))$, $\pi$ is onto. Hence, $({((E_7)^{\kappa, \mu})}^\gamma)_{\dot{F_1}(h{e_4})}/\Z_2 \cong SO(8)$. Therefore, ${((E_7)^{\kappa, \mu})}^\gamma)_{\dot{F_1}(h{e_4})}$ is isomorphism to $Spin(8)$ as a double covering group of $SO(8)$.
\vspace{3mm}

We shall determine the group structre of ${((E_7)^{\kappa, \mu})}^\gamma$.
\vspace{3mm}

{\sc Lemma 3.18.}{\it The Lie algebra ${((\gge_7)^{\kappa, \mu})}^\gamma$ of the group ${((E_7)^{\kappa, \mu})}^\gamma$ is given by }\\
$$
\begin{array}{l}
{((\gge_7)^{\kappa, \mu})}^\gamma
\vspace{1mm}\\
= \Big \{ {\mit\Phi}\Bigl(D_4  + D'_4 + \wti{A}_1(p) + i\pmatrix{\epsilon_1 & 0 & 0 \cr
                                                               0 & \epsilon_2 & q \cr
                                                               0 & \bar{q} & \epsilon_3}^{\sim}                                                                                                 , \pmatrix{0 & 0 & 0 \cr
                                                                           0 & \alpha_2 & x \cr
                                                                           0 & \ov{x} & \alpha_3 }
\vspace{1mm}\\
  -\tau \pmatrix{0 & 0 & 0 \cr                         
                 0 & \alpha_2 & x \cr                                                                                            0 & \ov{x} & \alpha_3 }, -\dfrac{3}{2}i{\epsilon_1} \Bigr) 
            \,\Big|\, D_4, D'_4 \in \so(4) \subset \so(8), \alpha_k \in C, p, q \in \H, 
\end{array}$$
$$
\begin{array}{l}
\qquad \qquad \qquad \qquad x \in \H^C, \epsilon_k \in \R, \epsilon_1+\epsilon_2+\epsilon_3 = 0 \Big \}.
 \end{array} $$   {\it In particular, we have} 
$$
  \dim({((\gge_7)^{\kappa, \mu})}^\gamma) = 34. $$
\vspace{-3mm}

{\sc Proposition 3.19.} ${((E_7)^{\kappa, \mu})}^\gamma \cong (Spin(4)\times Spin(8))/\Z_2,\,\,\Z_2 = \{(1, 1), (-1, \\
 -1) \}$.
\vspace{2mm}

{\sc Proof.} For $Spin(4) = Sp(1)\times Sp(1) = ({((E_7)^{\kappa, \mu})}^\gamma)_{\dot{F_1}(h),\ti{E}_1,\ti{E}_{-1}, \dot{E}_{23}}$(Proposi- \\
tions 3.1, 3.3) and $Spin(8) = ({((E_7)^{\kappa, \mu})}^\gamma)_{\dot{F_1}(h{e_4})}$\,(Proposition 3.17), we define a map $\phi_1 : Spin(4)\times Spin(8) \to {((E_7)^{\kappa, \mu})}^\gamma$ by 
$$
       \phi_1 (\alpha, \beta) = \alpha\beta.  $$
Then, $\phi_1$ is well-defined. For ${\mit\Phi}_4 \in \spin(4)$\,(Lemma 3.2) and ${\mit\Phi}_8 \in \spin(8)$\,(Lemma 3.14), since [${\mit\Phi}_4, {\mit\Phi}_8$] $= 0$, we have $\alpha\beta = \beta\alpha$. Hence, $\phi_1$ is a homomorphism. It is not difficult to see that $\Ker\phi_1 = \{(1, 1), (-1, -1)\} = \Z_2$. Since ${((E_7)^{\kappa, \mu})}^\gamma (\cong (Spin(12))^{\gamma}$ (see [4],[5])) is connected and $\dim({((E_7)^{\kappa, \mu})}^\gamma) = 34\,$(Lemma 3.18)$ = 6 + 28 = \dim(\spin(4)\oplus \spin(8)), \phi_1$ is onto. Thus, we have the required isomorphism $(Spin(4)\times Spin(8))/\Z_2 \cong {((E_7)^{\kappa, \mu})}^\gamma$.
\vspace{2mm}

Now, we shall determine the group structure of $(E_7)^{\sigma, \gamma}$.
\vspace{2mm}

{\sc Lemma 3.20.}{\it The Lie algebra $(\gge_7)^{\sigma, \gamma}$ of the group $(E_7)^{\sigma, \gamma}$ is given by }
$$
\begin{array}{l}
(\gge_7)^{\sigma, \gamma}
\vspace{1mm}\\
= \Big \{ {\mit\Phi}\Bigl(D_4 + D'_4 + \wti{A}_1(p) + i\pmatrix{\epsilon_1 & 0 & 0 \cr
                                                               0 & \epsilon_2 & q \cr
                                                               0 & \bar{q} & \epsilon_3}^{\sim}                                                                                                 , \pmatrix{\alpha_1 & 0 & 0 \cr
                                                                           0 & \alpha_2 & x \cr
                                                                           0 & \ov{x} & \alpha_3},
\end{array} $$
$$
\begin{array}{l}
 -\tau \pmatrix{\alpha_1 & 0 & 0 \cr                         
                 0 & \alpha_2 & x \cr                                                                                            0 & \ov{x} & \alpha_3}, \nu \Bigr) 
            \,\Big|\, D_4, D'_4 \in \so(4) \subset \so(8), \alpha_k \in C, p, q \in \H, x \in \H^C,
\vspace{1mm}\\ 
\qquad \qquad \qquad \qquad \qquad       \epsilon_k \in \R, \epsilon_1+\epsilon_2+\epsilon_3 = 0, \nu \in i\R \Big \}.  
\end{array}$$   {\it In particular, we have} 
$$
  \dim((\gge_7)^{\sigma, \gamma}) = 37. $$
\vspace{-3mm}

{\sc Proposition 3.21.} {\it For $A \in SU(2) = \{ A \in M(2, C) \, | \, (\tau{}^t\!A)A = E, \det A = 1 \}$, we define $C$-linear transformations $\phi(A)$ of $\gP^C$ by}
$$
\begin{array}{l}
      \phi(A)\Big(\pmatrix{\xi_1 & x_3 & \ov{x}_2 \cr
                             \ov{x}_3 & \xi_2 & x_1 \cr 
                             x_2 & \ov{x}_1 & \xi_3},
                    \pmatrix{\eta_1 & y_3 & \ov{y}_2 \cr
                             \ov{y}_3 & \eta_2 & y_1 \cr 
                             y_2 & \ov{y}_1 & \eta_3}, \xi, \eta \Big)
\end{array}$$
$$
\begin{array}{l}
\qquad \quad
     = \Big(\pmatrix{{\xi'_1} & {x'_3} & {\ov{x}_2}' \cr
                        {\ov{x}_3}' & {\xi'_2} & {x'_1} \cr
                         {x'_2} & {\ov{x}_1}' & {\xi'_3}},
               \pmatrix{{\eta'_1} & {y'_3} & {\ov{y}_2}' \cr
                        {\ov{y}_3}' & {\eta'_2} & {y'_1} \cr
                         {y'_2} & {\ov{y}_1}' & {\eta'_3}},
                          \xi', \eta'\Big),
\end{array}$$
\vspace{0.2mm}

$$
\begin{array}{c}
     \pmatrix{{\xi_1}' \cr \eta'} = A\pmatrix{\xi_1 \cr \eta}, \;\;
     \pmatrix{\xi' \cr {\eta'_1}} = A\pmatrix{\xi \cr \eta_1}, \;\;
     \pmatrix{{\eta'_ 2} \cr {\xi'_3}} = A\pmatrix{\eta_ 2 \cr \xi_3},
\end{array}$$
$$
\begin{array}{c}
     \pmatrix{{\eta'_3} \cr {\xi'_2}} = A\pmatrix{\eta_3 \cr \xi_2},\;\;
     \pmatrix{{x'_1} \cr {y'_1}} = (\tau A)\pmatrix{x_1 \cr y_1}, 
\end{array} $$
$$
\begin{array}{c}
     \pmatrix{{x'_2} \cr {y'_2}} = \pmatrix{x_2 \cr y_2}, 
     \pmatrix{{x'_3} \cr {y'_3}} = \pmatrix{x_3 \cr y_3}.
\end{array}$$
 {\it Then}, $\phi(A) \in (E_7)^{\sigma,\gamma}$.
\vspace{2mm}
 
{\sc Proof.} Let ${\mit\Phi} = {\mit\Phi}(2\nu E_1 \vee E_1, aE_1, - \tau aE_1, \nu), a \in C, \nu \in i\R$. Then, ${\mit\Phi} \in (\gge_7)^{\sigma,\gamma}$(Lemma 3.20). 
Therefore, for $A = \exp \pmatrix{\nu & a \cr
                                -{\tau a} & -\nu} \in SU(2)$, we have 
$\phi(A) = \exp{\mit\Phi} \in (E_7)^{\sigma,\gamma}$.
\vspace{3mm}

{\sc Proposition 3.22.} $(E_7)^\sigma \cong (SU(2) \times Spin(12))/\Z_2, \; \Z_2 = \{ (E, 1), (-E, $ $
 -\sigma) \} $.
\vspace{2mm}

{\sc Proof.}  The isomorphism is induced by the homomorphism $\varphi_1 : SU(2) \times Spin(12) \to (E_7)^{\sigma}$ by $\varphi_1(A, \delta) = \phi(A)\delta.$ (In detail, see [4], [5].)
\vspace{3mm}

{\sc Theorem 3.23.} $(E_7)^{\sigma, \gamma} \cong (SU(2) \times Spin(4) \times Spin(8))/(\Z_2 \times \Z_2),$ $
      \Z_2 \times \Z_2 = \{(E, 1, 1), (E, \sigma, \sigma) \} \times \{(E, 1, 1), (-E, \gamma, -\sigma\gamma) \}$.
\vspace{2mm}

{\sc Proof.}  For $SU(2)$\,(Proposition 3.21), $Spin(4) = ({((E_7)^{\kappa, \mu})}^\gamma)_{\dot{F_1}(h),\ti{E}_1,\ti{E}_{-1},\dot{E}_{23}}$ \\ 
(Propositions 3.1, 3.3) and $Spin(8) = ({((E_7)^{\kappa, \mu})}^\gamma)_{\dot{F_1}(h{e_4})}$(Proposition 3.17), we define a map $\varphi : SU(2)\times Spin(4)\times Spin(8) \to (E_7)^{\sigma, \gamma}$ by
$$
            \varphi(A, \alpha, \beta) = \phi(A)\alpha\beta. $$
Then, $\varphi$ is well-defined. From Propositions 3.19, 3.22, $\varphi$ is a homomorphim. We shall show that $\varphi$ is onto. Let $\rho \in (E_7)^{\sigma, \gamma}$. Since $(E_7)^{\sigma, \gamma} \subset (E_7)^\sigma$, there exist $A \in SU(2)$ and $\delta \in Spin(12)$ such that $\rho = \varphi_1 (A, \delta)$\,(Proposition 3.22). Now, From $\gamma\rho\gamma = \rho$, we have $\phi(A)(\gamma\delta\gamma) = \phi(A)\delta$. Hence,
$$
   \left\{\begin{array}{l}
             A = A
\vspace{1mm}\\
          \gamma\delta\gamma = \delta   \end{array} \right.
\quad \mbox{or} \quad
    \left\{\begin{array}{l}
             A = -A
\vspace{1mm}\\
          \gamma\delta\gamma = -\sigma\delta   \end{array} \right.  $$
The latter case is impossible because $A = 0$ is false. In the former case, from Proposition 3.19, there exist $\alpha \in Spin(4)$ and $\beta \in Spin(8)$ such that $\delta =  \phi_1(\alpha, \beta)$. Hence, we have 
\begin{eqnarray*}
    \rho \!\!\!& = &\!\!\!\varphi_1(A, \delta) = \phi(A)\delta = \phi(A)\phi_1(\alpha, \beta)
\vspace{1mm}\\
         \!\!\!& = &\!\!\!\phi(A)\alpha\beta = \varphi(A, \alpha, \beta).
\end{eqnarray*}
It is not difficult to see that
\begin{eqnarray*}
  \Ker \varphi \!\!\! & = & \!\!\! \{(E, 1, 1), (E, \sigma, \sigma), (-E, \gamma, -\sigma\gamma), (-E, \sigma\gamma, -\gamma) \}
\vspace{1mm}\\
              \!\!\! & = & \!\!\! \{(E, 1, 1),(E, \sigma, \sigma) \} \times \{(E, 1, 1), (-E, \gamma, -\sigma\gamma) \}
\vspace{1mm}\\
              \!\!\! & = & \!\!\! \Z_2 \times \Z_2.
\end{eqnarray*}
Thus, we have the required isomorphism $(SU(2) \times Spin(4) \times Spin(8))/(\Z_2 \times \Z_2) \cong (E_7)^{\sigma, \gamma }$.

\bigskip
\begin{flushright}
\begin{tabular}{l}
Toshikazu Miyashita \\
Tohbu High School \\
Agata, Tohbu, 389-0517, Japan \\
E-mail: serano@janis.or.jp \\
\end{tabular}
\end{flushright}

\end{document}